\documentclass[journal]{IEEEtran}
\usepackage{graphicx}
\usepackage{amsmath}
\usepackage{amssymb}
\usepackage{cite}

\usepackage[vlined,ruled]{algorithm2e}

\newtheorem{theorem}{Theorem}[section]

\newtheorem{corollary}[theorem]{Corollary}
\newtheorem{proposition}[theorem]{Proposition}

\newtheorem{remark}{Remark}

\newenvironment{pfof}[1]{\vspace{1ex}\noindent{\itshape Proof of
    #1:}\hspace{0.5em}} {\hfill\oprocend\vspace{1ex}}

\newcommand{\until}[1]{\{1,\dots, #1\}}

\newcommand{\setdef}[2]{\{#1 \; | \; #2\}}

\newcommand{\map}[3]{#1: #2 \rightarrow #3}

\newcommand{\real}{\mathbb{R}}

\newcommand\oprocendsymbol{\hbox{$\square$}}
\newcommand\oprocend{\relax\ifmmode\else\unskip\hfill\fi\oprocendsymbol}

\newcommand{\bsin}{\boldsymbol{\sin}}
\newcommand{\barcsin}{\boldsymbol{\arcsin}}
\newcommand{\bcos}{\boldsymbol{\cos}}

\renewcommand{\j}{\boldsymbol{\mathrm{j}}}

\newcommand{\vect}[1]{\mathbbold{#1}}
\newcommand{\vones}[1][]{\vect{1}_{#1}}
\newcommand{\vzeros}[1][]{\vect{0}_{#1}}
\DeclareSymbolFont{bbold}{U}{bbold}{m}{n}
\DeclareSymbolFontAlphabet{\mathbbold}{bbold}

\usepackage{enumerate}
\usepackage{xcolor}

\usepackage{booktabs}

\newcommand{\tb}{}

\definecolor{jg}{RGB}{14,120,14}

\graphicspath{{./Figures/}}

\ifCLASSINFOpdf
\else
\fi

\graphicspath{{./Figures/}}

\begin{document}
%
\title{Lossy DC Power Flow}
%
%
%

\author{John~W.~Simpson-Porco,~\IEEEmembership{Member,~IEEE}
\thanks{J.~W.~Simpson-Porco is with the Department of Electrical and Computer Engineering, University of Waterloo, Waterloo ON, N2L 3G1 Canada. Correspondance: {\texttt{jwsimpson@uwaterloo.ca}}.}
}

%
%

\markboth{Submitted to IEEE Transactions on Power Systems. This version: \today}%
{Shell \MakeLowercase{\textit{et al.}}: Bare Demo of IEEEtran.cls for Journals}
%



\maketitle

\begin{abstract}
The DC Power Flow approximation has been widely used for decades in both industry and academia due to its computational speed and simplicity, but suffers from inaccuracy, in part due to the assumption of a lossless network. Here we present a natural extension of the DC Power Flow to lossy networks. Our approach is based on reformulating the lossy active power flow equations into a novel fixed-point equation, and iterating this fixed-point mapping to generate a sequence of improving estimates for the active power flow solution. Each iteration requires the solution of a standard DC Power Flow problem with a modified vector of power injections. The first iteration returns the standard DC Power Flow, and one or two additional iterations yields a one or two order-of-magnitude improvement in accuracy. For radial networks, we give explicit conditions on the power flow data which guarantee (i) that the active power flow equations possess a unique solution, and (ii) that our iteration converges exponentially and monotonically to this solution. For meshed networks, we extensively test our results via standard power flow cases.
\end{abstract}



\section{Introduction}
\label{Sec:Introduction}



The power flow equations describe the steady-state transmission of power through an AC power grid, and are the basis for all power system analysis, operations, and control. These nonlinear equations are often included as equality constraints in power system optimization problems, or are solved as a subroutine within larger algorithms. While many problems only require the power flow equations to be solved once, applications such as contingency analysis and security-constrained dispatch/unit commitment require the repeated solution of many large power flow problems. 
This spectrum of requirements is matched by a spectrum of power flow models and solution algorithms, which range from exact to approximate. At one end of this spectrum are solution techniques such as Newton-Raphson and its derivatives. These algorithms precisely solve the power flow equations with no approximations, but are costly computationally, and their convergence is difficult to theoretically characterize \cite{JST-SAN:89,DM-DKM-KT:16}.

At the other end of this spectrum is the DC Power Flow (DCPF) approximation, which is an extensively used linearized model of active power flow; see \cite[Chapter 4.1.4]{AJW-BFW:96} for an introduction, and \cite{BS-JJ-OA:09} for a recent overview. The DCPF is approximate in several ways, but is also (i) linear (ii) easy to understand, (iii) computationally inexpensive, and (iv) guaranteed to converge, as it requires only the solution of a single set of sparse linear equations. Linearity has been extensively exploited for characterizing location marginal prices in market applications \cite{TJO-XC-YS:04}, and the DCPF is particularly attractive for combinatorial contingency and cascading failure studies \cite{VBCCDHPMZ:12}.

One of the main sources of inaccuracy in the DCPF model is that line resistances are neglected. The accuracy of this approximation was studied in \cite{KP-LM-DVD-RB:05}, and when compared to the exact power flow solution, the DCPF was shown to generate significant errors when $R/X$ ratios reach $1/4$. Similarly, \cite{BS-JJ-OA:09} observed inaccuracy when parallel transfer paths had significantly different $R/X$ ratios, and \cite{CC-PVH-RB:12} suggested that resistive line losses are the largest source of error in the DCPF. These errors also tend to be exacerbated in large networks \cite{TJO-XC-YS:04,BS-JJ-OA:09}.\footnote{Similarly, high $R/X$ ratios degrade the convergence of fast decoupled load flow; see \cite{DR-AB:88,MHH:93} and the references therein.}

As generalizations of the DCPF, various linearized models of power flow have recently been proposed in several different coordinate systems \cite{MEB-FFW:89b,TK-RC-MK:13,MF-CL-SL:13,PS-SB-MC:14,CC-PVH:14,DD-SB-MC:15,SB-FD:15,SB-SZ:15,SVD-SSG-YCC:15}. These models extend the DCPF to account for variations in bus voltages, but resistive losses can never be precisely accounted for within a linear/affine model, as they are inherently nonlinear. Attempts have nonetheless been made to heuristically incorporate losses into DCPF-type models, as is done in the so-called $\alpha$-matching model \cite{BS-JJ-OA:09}. This model begins from a known power flow solution, and treats the losses calculated from that known solution as an approximation to the losses for the DCPF solution \cite{BS-JJ-OA:09,YQ-DS-DT:12}; see also \cite{CC-PVH-RB:12} for a linear programming approach. While such heuristics can be effective, our goal is to formulate a more theoretically rigorous procedure for including resistive losses within a DCPF-inspired framework.

\smallskip


\subsection{Contributions}

This paper presents, analyzes, and tests the Lossy DC Power Flow (L-DCPF), a generalization of the classic DCPF to lossy networks. There are three main contributions. First, we derive the Lossy DCPF algorithm. Our approach is to reformulate the active power flow equations into a fixed-point form, and use this reformulation to define an iteration.
The first step of this iteration returns a DCPF-like approximation, while subsequent iterates yield improving approximations to the desired solution of the active power flow equations including resistive losses.
The method can be interpreted as defining a sequence of DC Power Flow problems, wherein the vector of power injections is updated at each step to account for losses, using the estimates of the phase angles from the previous iteration. For meshed networks, an additional correction term is required to enforce Kirchhoff's voltage law.

Second, we present a detailed theoretical analysis of the Lossy DCPF; this analysis is the major distinguishing feature of this work. We begin by establishing that if the Lossy DCPF iteration converges, it must converge to an active power flow solution. 
For radial (i.e., tree) networks, we then give explicit conditions on the problem data for convergence. 
Under these conditions, the iteration converges monotonically and exponentially to the unique solution of the active power flow equations. Moreover, the contraction mapping theorem explicitly provides (i) a bound on the distance from the true solution at every step, and (ii) a large set of initializations from which the iteration converges. These conditions can also be interpreted as sufficient conditions for the existence of a unique ``small-angle-difference'' power flow solution.\footnote{Casting power flow equations as fixed-point equations and applying contraction mapping principles is a classic approach \cite{FFW:77}, but has received renewed attention recently. See \cite{JWSP:17b} for a survey of results; for distribution systems, see, for example \cite{SB-SZ:15,MB-NG:16b}.}

Third and finally, we extensively test the Lossy DCPF on benchmark systems, ranging in size from 39 to 13,659 buses. For both base-case and highly-loaded systems, even one iteration of the L-DCPF provides an order-of-magnitude improvement in accuracy over the DCPF. {\tb The computational complexity and convergence rate of the L-MDCPF iteration is comparable to a chord Newton method, and the convergence is very robust with respect to initialization errors.}

\smallskip

The Lossy DCPF retains some aspects of the classic DCPF, but not others. In terms of similarities, it is simple, intuitive, requires minimal network data, and does not involve any state-dependent Jacobian matrices; one simply solves a sequence of DCPF problems with the same, constant, state-independent $B^\prime$ matrix. 
This approach can rigorously accommodate stressed networks and high $R/X$ ratios; the price we pay for these improvements is that the method is iterative, and each iteration requires the solution of a DCPF-like linear equation. 
At a conceptual level, this article builds a formal mathematical bridge \cite{NASEM:16} between the DC Power Flow approximation \cite{BS-JJ-OA:09} and the effects of resistive losses.
The theoretical analysis presented also extends and generalizes known results for lossless active power flow solvability, presented in \cite{FD-MC-FB:11v-pnas,FD-FB:13c}.



\subsection{Paper Organization} The remainder of this section establishes some notation. Section \ref{Sec:GridModel} clearly states the network model under consideration, and reviews the DCPF approximation along with a recently introduced modified DCPF. Section \ref{Sec:LDF-PF} introduces our Lossy DCPF model, first for a two-bus system in Section \ref{Sec:LDF-PF1} and then for general networks in Section \ref{Sec:LDF-PF2}. Theoretical convergence results for the L-DCPF iteration are detailed in Section \ref{Sec:LDF-PF3}. Section \ref{Sec:Numerics} provides some remarks on implementation before proceeding to extensive numerical testing. Section \ref{Sec:Conclusions} concludes and points to future research.

\subsection{Preliminaries and Notation}
\label{Sec:Notation}

\emph{Vectors, Matrices, Functions:} We let $\real$ denote the set of real numbers. Given $x \in \real^{n}$, $[x] \in \real^{n\times n}$ is the diagonal matrix with $x$ on the diagonal, and $\|x\|_{\infty} = \max_{i}|x_i|$. 
Throughout, $\vones[n]$ and $\vzeros[n]$ are the $n$-dimensional vectors of unit and zero entries. The $n \times n$ identity matrix is $I_n$. For a matrix $M \in \real^{n \times m}$, $M^{\sf T}$ is its transpose, $\mathrm{im}(M)$ is its  range, $\mathrm{ker}(M)$ is its null-space, and $\|M\|_{\infty}$ is its induced $\infty$-norm. For $x \in \real^n$, we define the vector functions $\bsin(x) \triangleq (\sin(x_1),\ldots,\sin(x_n))^{\sf T}$, $\barcsin(x)$ as the inverse of $\bsin(x)$, $\bcos(x) \triangleq (\cos(x_1),\ldots,\cos(x_n))^{\sf T}$, $x^2 = (x_1^2,\ldots,x_n^2)^{\sf T}$, and for $x$ with nonnegative components, $\sqrt{x} = (\sqrt{x_1},\ldots,\sqrt{x_n})^{\sf T}$.

\smallskip

\emph{Graphs and graph matrices:} A graph is a pair $(\mathcal{N},\mathcal{E})$, where $\mathcal{N}$ is the set of nodes and $\mathcal{E} \subseteq \mathcal{N} \times \mathcal{N}$ is the set of edges; we consider only connected graphs. 
If a label $e \in \until{|\mathcal{E}|}$ and an arbitrary orientation is assigned to each edge $e \sim (i,j) \in \mathcal{E}$, the \emph{node-edge incidence matrix} $A \in \real^{|\mathcal{N}|\times|\mathcal{E}|}$ is defined component-wise as $A_{ke} = 1$ if node $k$ is the source node of edge $e$ and as $A_{ke} = -1$ if node $k$ is the sink node of edge $e$, with all other elements being zero. For $x \in \real^{|\mathcal{N}|}$, $A^{\sf T}x \in \real^{|\mathcal{E}|}$ is the vector with components $x_i-x_j$, with $(i,j) \in \mathcal{E}$. If the graph is radial (a tree), then $\mathrm{ker}(A) = \emptyset$.
The space $\mathrm{ker}(A)$ is the cycle-space of the graph, and we let $c$ be its dimension. We let $C \in \{-1,0,1\}^{|\mathcal{E}|\times c}$ denote the \emph{(oriented) edge-cycle incidence matrix} \cite[Section 3]{TK-CL-KM-DM-RR-TU-KAZ:09}; the columns of $C$ index cycles of the network, and span the cycle space, i.e., $AC = \vzeros[]$.
The matrix $|A|$ \textemdash{} the ``unsigned'' incidence matrix \textemdash{} is obtained by setting all non-zero elements of $A$ equal to $1$.
If a single row of $A$ is removed, we call the result a reduced incidence matrix $A_r$, which satisfies $\mathrm{ker}(A_r^{\sf T}) = \emptyset$. 
%
%
%
If a weight $w_{ij} > 0$ is assigned to each edge $(i,j)\in\mathcal{E}$, then ${\sf L} = {\sf L}^{\sf T} = A_r[w]A_r^{\sf T}$ is the corresponding \emph{reduced (or ``grounded'') Laplacian matrix}, and is positive definite.

\section{Power Flow Equations, DC Power Flow, and Modified DC Power Flow}
\label{Sec:GridModel}

\subsection{Grid Model and Power Flow Equations}

We consider the positive sequence representation of a balanced three-phase network, modeled as a weighted graph $(\mathcal{N},\mathcal{E})$ where $\mathcal{N} = \{1,\ldots,n+1\}$ is the set of $n+1$ buses and $\mathcal{E} \subset \mathcal{N} \times \mathcal{N}$ is the set of $m$ branches. The admittance matrix of the network is denoted by $Y = G + \j B$, where $G$ is the conductance matrix and $B$ is the susceptance matrix. The elements of $Y$ are $Y_{ij} = -y_{ij}/t_{ij}$, where $y_{ij}$ is the transfer admittance between buses $i$ and $j$ \textemdash{} including off-nominal tap ratios $t_{ij} \geq 1$ \textemdash{} while $Y_{ii}$ is defined as $Y_{ii} = -\sum_{j=1}^{n+1} y_{ij}/t_{ij}^2$. To every bus $i$ we associate a phasor voltage $V_i \angle \theta_i$ and a complex power $P_i + \j Q_i$, which are related to each other through the power flow equations
\begin{subequations}
\begin{align}
\label{Eq:Active}
P_i &= \sum_{j=1}^{n+1} V_iV_j\left(G_{ij}\cos(\theta_i-\theta_j)+B_{ij}\sin(\theta_i-\theta_j)\right)\\
\label{Eq:Reactive}
Q_i &= \sum_{j=1}^{n+1} V_iV_j\left(G_{ij}\sin(\theta_i-\theta_j)-B_{ij}\cos(\theta_i-\theta_j)\right)
\end{align}
\end{subequations}
for each $i \in \{1,\ldots,n+1\}$. We will take bus $n+1$ as the slack bus with $\theta_{n+1} = 0$, $V_{n+1}$ fixed, and $P_{n+1}, Q_{n+1}$ unknown. In this paper we focus on \eqref{Eq:Active} and not on \eqref{Eq:Reactive}. Buses $\{1,\ldots,n\}$ are considered as $\mathsf{PV}$ buses at which $P_i$ and $V_i$ are fixed and $\theta_i$ is unknown. Since voltage magnitudes are assumed fixed at buses $\until{n}$, constant power and ZIP load models are equivalent, so we may \textemdash{} without loss of generality \textemdash{} absorb any shunt-conductance/constant-current active power load models into the constant power $P_i$ in \eqref{Eq:Active}.

\subsection{DC and Modified DC Power Flow}
\label{Sec:DC-and-MDC}

The model \eqref{Eq:Active}--\eqref{Eq:Reactive} is often simplified further in a number of ways. Under the assumptions that (i) all conductances are zero (ii) voltage magnitudes are known and constant, and (iii) reactive power is unimportant, the model  \eqref{Eq:Active}--\eqref{Eq:Reactive} may be simplified to the lossless decoupled model
\begin{equation*}
P_i = \sum_{j=1}^{n+1}\nolimits V_iV_jB_{ij}\sin(\theta_i-\theta_j)\,,\qquad i \in \until{n}\,.
\end{equation*}
We write this lossless model in vector notation as
\begin{equation}\label{Eq:LosslessActive}
P_r = A_r{\sf D}_{B}\,\bsin(A_r^{\sf T}\theta_r)\,,
\end{equation}
where $P_r = (P_1,\ldots,P_n)^{\sf T}$ and $\theta_r = (\theta_1,\ldots,\theta_n)^{\sf T}$ are column vectors of power injections and angles for buses $\{1,\ldots,n\}$, $A_r$ is the reduced bus-branch incidence matrix of the network obtained by removing the row of $A$ corresponding to the slack bus, and
\begin{equation}\label{Eq:DB}
{\sf D}_B \triangleq [V_iV_jB_{ij}]
\end{equation}
is a diagonal and positive-definite $m \times m$  matrix of branch weights. 
{
The \emph{DC Power Flow} (DCPF) assumes\footnote{The terminology ``DC Power Flow'' is usually packaged with several additional assumptions. Under these assumptions, ${\sf D}_{B}$ is simply the diagonal matrix of inverse series line reactances, and the common notation for ${\sf L}_{B}$ is then $B^\prime$. Although we do not impose these additional assumptions, in this article we nonetheless refer to \eqref{Eq:DCPF} as the DCPF.} that phase angle differences $A_r^{\sf T}\theta_r$ are small, and hence that $\bsin(A_r^{\sf T}\theta_r) \approx A_r^{\sf T}\theta_r$. 
The power flow equations \eqref{Eq:LosslessActive} therefore reduce to a system of linear equations $P_r \approx {\sf L}_B\theta_r$, where
\begin{equation}\label{Eq:Lb}
{\sf L}_{B} \triangleq A_r{\sf D}_BA_r^{\sf T}
\end{equation}
is the reduced Laplacian matrix associated with the branch weights in ${\sf D}_B$. We therefore obtain the explicit DCPF solution
\begin{equation}\label{Eq:DCPF}
\theta_{r,\mathrm{DC}} = {\sf L}_B^{-1}P_r\,.
\end{equation}
The \emph{Modified DC Power Flow} (MDCPF), proposed in \cite{FD-FB:13c}, gives an improvement in accuracy and rigor over the DCPF. The idea is to make the substitution $\psi = \bsin(A_r^{\sf T}\theta_r)$ in \eqref{Eq:LosslessActive}, leading to the pair of equations
\begin{subequations}\label{Eq:MDCPF-Unsolved}
\begin{align}
\label{Eq:MDCPF-Unsolved-a}
P_r &= A_r{\sf D}_{B}\psi\\
\label{Eq:MDCPF-Unsolved-b}
\psi &= \bsin(A_r^{\sf T}\theta_r)
\end{align}
\end{subequations}
The general solution to \eqref{Eq:MDCPF-Unsolved-a} is
\begin{equation}\label{Eq:Psi0}
\psi = \underbrace{A_r^{\sf T}{\sf L}_{B}^{-1}P_r}_{:= \psi_{\rm MDC}} + {\sf D}_{B}^{-1}Cx\,,
\end{equation}
where $C$ is the edge-cycle incidence matrix and $x \in \real^{c}$. The first term $\psi_{\rm MDC}$ is a particular solution, while the second term ${\sf D}_{B}^{-1}Cx$ parameterizes the homogeneous solution. This is easily verified, since
$$
A_r{\sf D}_B\psi_{\rm MDC} = A_r{\sf D}_BA_r^{\sf T}{\sf L}_B^{-1}P_r = {\sf L}_B{\sf L}_B^{-1}P_r = P_r\,,
$$
and $A_r{\sf D}_{B}({\sf D}_B^{-1}Cx) = A_rCx = \vzeros[n]$. Assuming that $x = \vzeros[c]$ and $\|\psi_{\rm MDC}\|_{\infty} < 1$, one may apply $\barcsin$ to both sides of \eqref{Eq:MDCPF-Unsolved-b} to obtain the solution $A_r^{\sf T}\theta_r = \barcsin(\psi_{\rm MDC})$, or
\begin{equation}\label{Eq:MDCPF}
\begin{aligned}
\theta_{r,\mathrm{MDC}} &= (A_rA_r^{\sf T})^{-1}A_r\,\barcsin(A_r^{\sf T}{\sf L}_{B}^{-1}P_r)\,.
\end{aligned}
\end{equation}
The modified solution \eqref{Eq:MDCPF} is more accurate than the DCPF in stressed networks, due to the explicit use of the nonlinearity $\barcsin(\cdot)$, and \eqref{Eq:MDCPF} is provably an exact solution for radial networks \cite{FD-FB:13c}.
Finally, note that if one approximates $\barcsin(x) \approx x$, then \eqref{Eq:MDCPF} reduces to the DCPF \eqref{Eq:DCPF}. A more a detailed numerical comparison of the DCPF and the MDCPF may be found in \cite{FD-FB:13c}. 
Our goal now is to extend the classic DCPF \eqref{Eq:DCPF} and the Modified DCPF \eqref{Eq:MDCPF} to lossy networks.
}

\section{Lossy DC Power Flow}
\label{Sec:LDF-PF}

This section presents a derivation of the (Modified) Lossy DCPF. We begin with the two-bus case in Section \ref{Sec:LDF-PF1}, before proceeding to the general case in Section \ref{Sec:LDF-PF2}.

\subsection{Derivation for Two-Bus Case}
\label{Sec:LDF-PF1}

Consider the two-bus network consisting of two \textsf{PV} buses connected by a branch of admittance $y = g - \boldsymbol{\mathrm{j}}b$. With bus 2 as the slack bus, the active power flow equations \eqref{Eq:Active} are
\begin{subequations}\label{Eq:TwoBusAPF}
\begin{align}\label{Eq:TwoBusAPFa}
P_1 &= gV_1^2 - gV_1V_2\cos\theta + bV_1V_2\sin\theta\\
\label{Eq:TwoBusAPFb}
P_2 &= gV_2^2 - gV_1V_2\cos\theta - bV_1V_2\sin\theta\,,
\end{align}
\end{subequations}
where $g, b > 0$ and $\theta = \theta_1$. Making the substitution $\psi = \sin(\theta)$, we rewrite \eqref{Eq:TwoBusAPFa} as
\begin{subequations}\label{Eq:TwoBusPsiForm}
\begin{align}\label{Eq:TwoBusPsiForma}
P_1 &= gV_1^2 - gV_1V_2\sqrt{1-\psi^2} + bV_1V_2\psi\\
\label{Eq:TwoBusPsiFormb}
\psi &= \sin(\theta)\,,
\end{align}
\end{subequations}
%
%
%
%
%
%
where we have assumed that $|\theta| < \frac{\pi}{2}$ and used that $\cos \theta = (1-\sin^2\theta)^{1/2} = (1-\psi^2)^{1/2}$. Our goal is now to solve \eqref{Eq:TwoBusPsiForma} for $\psi$, since we can then recover the phase angle as $\theta = \arcsin(\psi)$. We divide \eqref{Eq:TwoBusPsiForma} by $bV_1V_2$ and rearrange to obtain
\begin{equation}\label{Eq:PsiFP-TwoBus}
\psi = \frac{P_1}{bV_1V_2} + \frac{g}{b}\left(\sqrt{1-\psi^2} - \frac{V_1}{V_2}\right)\,.
\end{equation}
We will think of \eqref{Eq:PsiFP-TwoBus} as defining an update equation for $\psi$. One substitutes an estimate  of $\psi$ into the right-hand side of \eqref{Eq:PsiFP-TwoBus}, and the equation returns an updated estimate. Since $\psi = \sin\theta$, our initial estimate would be the ``flat start'' condition $\psi[0] = \sin(0) = 0$. Substituting this into \eqref{Eq:PsiFP-TwoBus} generates
$$
\psi[1] = P_1/(bV_1V_2) + (g/b)(1-V_1/V_2)\,,
$$
which equals $\psi_{\rm MDC} = P_1/(bV_1V_2)$ when $V_1 = V_2$. Therefore, the first iteration $\psi[1]$ returns a variation of the Modified DC Power Flow.
Repeating this process, one can generate a sequence of approximations. We define the \emph{Lossy Modified DC Power Flow} (L-MDCPF) iteration as
\begin{equation}\label{Eq:MLDCPF-TWO}
\boxed{
\begin{aligned}
\psi[k] &= \frac{P_1}{bV_1V_2} + \frac{g}{b}\left(\sqrt{1-(\psi[k-1])^2}-\frac{V_1}{V_2}\right)\,\\
{\theta}[k] &= \arcsin(\psi[k])\,.
\end{aligned}
}
\end{equation}
Once again, the terminology ``modified'' refers to phase angles being calculated using an $\arcsin$ \cite{FD-FB:13c}, in contrast to standard DC Power Flow. If we approximate $\arcsin(\psi) \approx \psi$ in \eqref{Eq:MLDCPF-TWO}, we instead obtain the \emph{Lossy DC Power Flow} (L-DCPF)
\begin{equation}\label{Eq:LDCPF-TWO}
\boxed{
\begin{aligned}
{\theta}[k] &= \frac{P_1}{bV_1V_2} + \frac{g}{b}\left(\sqrt{1-({\theta}[k-1])^2}-\frac{V_1}{V_2}\right)\,.
\end{aligned}
}
\end{equation}

%

\subsection{Derivation for General Meshed Networks}
\label{Sec:LDF-PF2}

We now extend the development of Section \ref{Sec:LDF-PF1} to arbitrary networks. As we will see, a complicating factor not present in the radial two-bus case is that Kirchhoff's voltage law must be satisfied around every cycle of the network. We begin by returning to the active power flow \eqref{Eq:Active}, which we write as
\begin{equation}\label{Eq:LossyActiveWrittenOut}
\begin{aligned}
P_i &=  \sum_{j\neq i} \nolimits B_{ij}V_iV_j\sin(\theta_i-\theta_j) + G_{ii}V_i^2 \\& \quad + \sum_{j\neq i} \nolimits G_{ij}V_iV_j\cos(\theta_i-\theta_j)\,,
\end{aligned}
\end{equation}
for $i \in \until{n}$. Just as we did in Section \ref{Sec:DC-and-MDC}, this may be written in vector form as
\begin{equation}\label{Eq:ActiveGeneralVector}
\begin{aligned}
P_r &= A_r\mathsf{D}_{B}\,\bsin(A_r^{\sf T}\theta_r) + G_{\rm diag}V_r^2\\ &\qquad - |A|_r\mathsf{D}_{G}\,\bcos(A_r^{\sf T}\theta_r)\,.
\end{aligned}
\end{equation}
where $G_{\rm diag} \in \real^{n \times n}$ is the diagonal matrix with elements $G_{ii}$, $V_r^2 = (V_1^2,\ldots,V_n^2)^{\sf T}$, $|A|_r$ is the ``unsigned'' reduced bus-branch incidence matrix (Section \ref{Sec:Notation}), and
\begin{equation}\label{Eq:DG}
\mathsf{D}_G = [V_iV_j|G_{ij}|]
\end{equation}
is an $m \times m$ diagonal matrix of branch weights.
We again introduce the change of variables $\psi \triangleq \bsin(A_r^{\sf T}\theta_r) \in \real^m$. Assuming that the branch-wise angle differences $A_r^{\sf T}\theta_r$ satisfy $\|A_r^{\sf T}\theta_r\|_{\infty} < \frac{\pi}{2}$, it follows from $\sin^2(z) + \cos^2(z) = 1$ that
$$
\bcos(A_r^{\sf T}\theta) = \sqrt{\vones[m]-\psi^2}\,,
$$
where $\psi^2 = (\psi_1^2,\ldots,\psi_m^2)^{\sf T}$. Eliminating $\theta_r$ from \eqref{Eq:ActiveGeneralVector}, we arrive at the equivalent pair of equations
\begin{subequations}\label{Eq:ActivePsi}
\begin{align}
P_r &= A_r\mathsf{D}_{B}\psi + G_{\rm diag}V_r^2 - |A|_r\mathsf{D}_{G}\sqrt{\vones[m]-\psi^2}\label{Eq:ActivePsia}\\
\label{Eq:ActivePsib}
\psi &= \bsin(A_r^{\sf T}\theta_r)\,.
\end{align}
\end{subequations}
If one temporarily thinks of the third term on the right-hand side of \eqref{Eq:ActivePsia} as being constant, then \eqref{Eq:ActivePsia} is linear in $\psi$. By mirroring our arguments from Section \ref{Sec:DC-and-MDC}, we can therefore write the general solution to \eqref{Eq:ActivePsia} as
\begin{equation}\label{Eq:PsiFP-General}
\begin{aligned}
\psi &= A_r^{\sf T}{\sf L}_B^{-1}\Big(P_r - G_{\rm diag}V_r^2 \\ & \qquad \qquad \qquad + |A|_r\mathsf{D}_{G}\sqrt{\vones[m]-\psi^2}\,\Big) + {\sf D}_B^{-1}Cx\,,
\end{aligned}
\end{equation}
where $C$ is the edge-cycle incidence matrix and $x \in \real^c$ is to be determined. Once $\psi$ is known, the phase angle differences $A_r^{\sf T}\theta_r = \barcsin(\psi)$ are calculated from \eqref{Eq:ActivePsib}. However, by Kirchhoff's voltage law, these differences must sum to zero ($\mathrm{mod}\,2\pi$) around every cycle in the network, and therefore
\begin{equation}\label{Eq:PsiFP-General2}
\begin{aligned}
C^{\sf T}\barcsin(\psi) = \vzeros[c]\,\, (\mathrm{\bf mod}\,2\pi)\,.
\end{aligned}
\end{equation}
To summarize so far, the active power flow equation \eqref{Eq:LossyActiveWrittenOut} is equivalent to the pair \eqref{Eq:ActivePsi}, which is in turn equivalent to the pair \eqref{Eq:PsiFP-General}--\eqref{Eq:PsiFP-General2}. We now think of \eqref{Eq:PsiFP-General} as defining an update equation for $\psi$. If we substitute the ``flat-start'' initialization $\psi[0] = \vzeros[m], x[0] = \vzeros[c]$ into \eqref{Eq:PsiFP-General}, we obtain the iterate
\begin{equation}\label{Eq:Psi111}
\psi[1] = \psi_{\rm MDC} + A_r^{\sf T}{\sf L}_B^{-1}(|A|_r\mathsf{D}_{G}\vones[m]-G_{\rm diag}V_r^2)\,.
\end{equation}
When $V_i = V_j$ for all buses $i,j$, the quantity within the brackets in $\psi[1]$ vanishes, and $\psi[1] = \psi_{\rm MDC}$. Therefore, just as in the two-bus case of Section \ref{Sec:LDF-PF1}, the first iteration returns a variation of the Modified DCPF solution. 
%
%
%
Repeating this process, one can generate a sequence of approximations. For networks with cycles, the variable $x$ must also be updated by solving \eqref{Eq:PsiFP-General2}. While many options are possible for this update, we propose the following chord Newton step
$$
x[k+1] = x[k] - (C^{\sf T}{\sf D}_B^{-1}C)^{-1}C^{\sf T}\barcsin(\psi[k])\,.
$$
In general then, we can define the \emph{Lossy Modified DC Power Flow} (L-MDCPF) iteration for meshed networks as
\begin{equation}\label{Eq:MLDCK-PF-General}
\boxed{
\begin{aligned}
x[k&+1] = x[k] - (C^{\sf T}{\sf D}_B^{-1}C)^{-1}C^{\sf T}\barcsin(\psi[k])\\
\psi[k&+1] =  A_r^{\sf T}{\sf L}_{B}^{-1}\bigg(P_r - G_{\rm diag}V_r^2\\  
&+ |A|_r\mathsf{D}_G\sqrt{\vones[m]-(\psi[k])^2}\,\bigg) + {\sf D}_{B}^{-1}Cx[k+1]\\
{\theta}_r[k&+1] = (A_rA_r^{\sf T})^{-1}A_r\,\barcsin(\psi[k+1])\,.
\end{aligned}
}
\end{equation} 
For radial networks, one simply removes the iteration for $x$ and sets $x = \vzeros[c]$. If we (i) assume $x = \vzeros[c]$, and (ii) approximate $\barcsin(\psi) \approx \psi$, then we may eliminate $\psi[k]$ and state the simpler \emph{Lossy DC Power Flow} (L-DCPF)
\begin{equation}\label{Eq:LDCK-PF-General}
\boxed{
\begin{aligned}
{\theta}_r[k+1] &\triangleq  {\sf L}_{B}^{-1}\bigg(P_r - G_{\rm diag}V_r^2\\  & + |A|_r\mathsf{D}_G\sqrt{\vones[m]-(A_r^{\sf T}{\theta}[k])^2}\,\bigg)\,.
\end{aligned}
}
\end{equation}

\begin{remark}[\bf Interpretation of Lossy DCPF]\label{Rem:Interp}
In \eqref{Eq:PsiFP-General}, one can interpret the quantity in brackets as an \emph{effective} vector of power injections, which includes the effects of resistive losses. Unfortunately, this effective injection vector depends on $\psi$, which is the variable we are trying to determine. The L-MDCPF and L-DCPF iterations obtain an estimate for this effective power injection vector by using the value of $\psi$ (or equivalently, the estimated phase angles) from the previous iteration (see also Section \ref{Sec:Computational}). This is an automatic loss-allocation procedure, which exploits the physics of power flow, in the spirit of an iterated $\alpha$-matching procedure \cite[Eqns. (6)-(7)]{YQ-DS-DT:12}. One can therefore think of the iterations \eqref{Eq:MLDCK-PF-General} and \eqref{Eq:LDCK-PF-General} as solving a sequence of DC Power Flows, where each iteration provides the hot-start point for the next. The additional variable $x[k]$ ensures that the solution $A_r^{\sf T}\theta_r = \barcsin(\psi)$ satisfies Kirchhoff's voltage law.
\oprocend
\end{remark}

\smallskip

\section{Theoretical Results: Power Flow Solutions and Convergence of Lossy DCPF}
\label{Sec:LDF-PF3}

The Lossy Modified DC Power Flow \eqref{Eq:MLDCK-PF-General} was formulated by transforming the auxiliary equation \eqref{Eq:ActivePsia} into \eqref{Eq:PsiFP-General}, then turning \eqref{Eq:PsiFP-General} into an iteration. If this iteration converges to a power flow solution, we can be confident that by taking $k$ large enough, ${\theta}_r[k]$ will be a good approximation to the power flow solution. Moreover, if this iteration converges \emph{monotonically}, we can be confident that each successive approximation is more accurate than the previous one. We now present results in this direction; all proofs are contained in Appendix \ref{App:1}.

\smallskip

\begin{proposition}[\bf Lossy Modified DCPF and Active Power Flow]\label{Thm:GeneralConvergence}
%
If the sequences $\{x[k]\}$ and $\{\psi[k]\}$ generated by the Lossy Modified DC Power Flow \eqref{Eq:MLDCK-PF-General} converge to values $x^* \in \real^c$ and $\psi^* \in \real^m$ respectively, then the sequence $\{A_r^{\sf T}\theta_r[k]\}$ of phase differences converges to a solution $A_r^{\sf T}\theta_r^* = \barcsin(\psi^*) \,\, (\mathrm{\bf mod}\,2\pi)$ of the active power flow equations \eqref{Eq:LossyActiveWrittenOut}.
%
\end{proposition}

\smallskip

Proposition \ref{Thm:GeneralConvergence} simply says that the L-MDCPF is an iterative algorithm for solving the active power flow equations \eqref{Eq:LossyActiveWrittenOut}. The iteration variable $x[k]$ is required due to the possibility of loop flows in meshed networks \cite{AJK:72}. When the network is radial, we may set $x = \vzeros[c]$ and the situation simplifies.


\smallskip

\begin{corollary}[\bf Lossy Modified DCPF and Active Power Flow on Radial Networks]\label{Cor:GeneralConvergence}
If the network is radial and the sequence $\{\psi[k]\}$ from the Lossy Modified DC Power Flow \eqref{Eq:MLDCK-PF-General} converges to a value $\psi^* \in \real^m$, then the sequence $\{\theta_r[k]\}$ converges to the solution $\theta_r^* = (A_rA_r^{\sf T})^{-1}A_r\barcsin(\psi^*)$ of the active power flow equations \eqref{Eq:LossyActiveWrittenOut}.
\end{corollary}

\smallskip

Proposition \ref{Thm:GeneralConvergence} and Corollary \ref{Cor:GeneralConvergence} address what happens \emph{if} the sequences of iterates $\{x[k],\psi[k]\}$ converge. We now seek conditions on the network data for which convergence is guaranteed. If we are to converge to a solution however, we should ensure there is a solution to converge to in the first place, so we now address this point. We restrict ourselves to radial networks, as the meshed case remains an open problem even for lossless networks \cite{FD-MC-FB:11v-pnas}. For simplicity, we assume all voltage magnitudes are equal and that there are no off-nominal tap ratios; these latter assumptions are easily relaxed.

\smallskip

\begin{theorem}[\bf Existence and Uniqueness of Solutions to Active Power Flow on Radial Networks]\label{Thm:Conditions}
Consider the lossy active power flow equations \eqref{Eq:LossyActiveWrittenOut}. Assume the network is radial and that $V_i = V > 0$ for all buses $i \in \until{n+1}$. Let
\begin{equation}\label{Eq:rho}
\rho \triangleq \|{\sf D}_{B}^{-1}A_r^{-1}|A|_r{\sf D}_G\|_{\infty}\\
\end{equation}
be a measure of the $R/X$ ratio of the network, let $\psi_{\rm MDC} = A_r^{\sf T}{\sf L}_B^{-1}P_r$ be the auxiliary solution to the Modified DC Power Flow from \eqref{Eq:Psi0}, and let
\begin{equation}\label{Eq:Gamma}
\Gamma \triangleq \|\psi_{\rm MDC}\|_{\infty}
\end{equation}
measure the size of that solution. If
\begin{equation}\label{Eq:Condition}
\Gamma^2 + 2\Gamma\rho < 1\,,
\end{equation}
then the active power flow equations \eqref{Eq:LossyActiveWrittenOut} possess a unique solution $\theta_r^* = (\theta_1^*,\ldots,\theta_n^*)$ with branch-wise phase differences $A_r^{\sf T}\theta_r^* \in \real^m$ satisfying $\|A_r^{\sf T}\theta_r^*\|_{\infty} \leq \arcsin(\beta_-) < \frac{\pi}{2}$, where
$$
\beta_{\pm} = \frac{\Gamma + \rho}{1+\rho^2} \pm \frac{\rho}{1+\rho^2}\sqrt{1-(\Gamma^2+2\Gamma\rho)}
$$
satisfy $0 \leq \beta_- < \beta_+ \leq 1$. Moreover, there \emph{do not} exist any solutions with phase differences satisfying
\begin{equation}\label{Eq:NoSolutions}
\arcsin(\beta_-) < |\theta_i-\theta_j| < \arcsin(\beta_+)\,,\
\end{equation}
for any branch $(i,j)\in\mathcal{E}$.
\end{theorem}

\smallskip

The main condition \eqref{Eq:Condition} of Theorem \ref{Thm:Conditions} is quite intuitive. It says that the Modified DC Power Flow solution $\psi_{\rm MDC}$ \textemdash{} and hence the phase angle differences $A_r^{\sf T}\theta_{r,\mathrm{MDC}} = \barcsin(\psi_{\rm MDC})$ that result from it \textemdash{} should not be too large, and that neither should the network $R/X$ ratio, as measured by $\rho$. These two variables play off one-another: for fixed $\Gamma \in (0,1)$, the condition \eqref{Eq:Condition} says that $\rho < \frac{1-\Gamma^2}{2\Gamma}$. This permits high $R/X$ ratios when the MDCPF solution has small angles, and lower $R/X$ ratios as the angles from the MDCPF solution grow.\footnote{A similar observation regarding $R/X$ ratios was made in \cite{FFW:77}.} Moreover, in the two-bus case of Section \ref{Sec:LDF-PF1}, \eqref{Eq:Condition} is both necessary and sufficient for the existence and uniqueness of a solution. When $\rho = 0$, Theorem \ref{Thm:Conditions} reduces to the lossless case studied in \cite{FD-MC-FB:11v-pnas}. The final equation \eqref{Eq:NoSolutions} implies that $\theta_r^*$ is the unique solution within a large region in angle-space.

\smallskip

In the proof of Theorem \ref{Thm:Conditions}, we show that the condition \eqref{Eq:Condition} implies that the right-hand side of \eqref{Eq:PsiFP-General} is a contraction mapping. Since the Lossy Modified DCPF \eqref{Eq:MLDCK-PF-General} is based on iterating this mapping, we have a strong convergence result.

\smallskip

\begin{corollary}[\bf Convergence of Lossy Modified DC Power Flow for Radial Networks]\label{Cor:RadialConvergence}
Consider the Lossy Modified DC Power Flow \eqref{Eq:MLDCK-PF-General}, and assume the conditions of Theorem \ref{Thm:Conditions} hold. Then from every initialization ${\psi}[0] \in (-\beta_+,\beta_+)^n$, the corresponding sequence of phase angle estimates $\{{\theta}_r[k]\}_{k=1,2,\ldots}$ converges monotonically and exponentially to the unique solution $\theta_r^*$ of the active power flow equations \eqref{Eq:LossyActiveWrittenOut} satisfying $\|A_r^{\sf T}\theta_r^*\|_{\infty} \leq \arcsin(\beta_-)$. The sequence of variables $\psi[k]$ satisfy the error estimate
$$
\|\psi[k] - \psi^*\|_{\infty} \leq \frac{\Gamma}{1-c}c^k\,,
$$
where $\psi^* = \bsin(A_r^{\sf T}\theta_r^*)$ and $c = \rho\beta_-(1-\beta_-^2)^{-\frac{1}{2}} < 1$.
\end{corollary}

\smallskip

Corollary \ref{Cor:RadialConvergence} delivers the basic intuition for why the Lossy DC Power Flow should work well: the iteration \eqref{Eq:MLDCK-PF-General} is, under reasonable conditions, a contraction mapping. Monotonic convergence tells us that each iteration improves on the previous one, and exponential convergence tells us this improvement is fairly rapid. The result provides bounds on both the convergence rate and the set of convergent initializations. Corollary \ref{Cor:RadialConvergence} is restricted to radial networks; in Section \ref{Sec:Numerics} we will test the Lossy DCPF iteration on meshed networks.

\section{Numerical Testing}
\label{Sec:Numerics}

Having established the theoretical basis for our Lossy DC Power Flow formulation, we now proceed to numerical testing. In Section \ref{Sec:Computational} we comment on how to properly compute the L-MDCPF, before proceeding to testing in Section \ref{Sec:Simulations}.

\subsection{Computational Considerations}
\label{Sec:Computational}

The Lossy Modified DC Power Flow \eqref{Eq:MLDCK-PF-General} is written out explicitly in terms of inverses. It is however computationally preferable to code the algorithm in terms of linear equations, so that sparse solution techniques can be applied. Algorithm \ref{Alg:LDCPF} details one reasonable implementation for the L-MDCPF \eqref{Eq:MLDCK-PF-General}; other implementations are obviously possible. The implementation of the simpler Lossy DCPF \eqref{Eq:LDCK-PF-General} is similar.

\begin{algorithm}

\SetKwInput{KwData}{Inputs}
\SetKwInput{KwResult}{Outputs}
\SetKwFunction{Solve}{Solve}

 \KwData{grid data, max iterations $I_{\rm max}$, $\psi[0]$, $x[0]$}
 \KwResult{voltage phase angles ${\theta}_r$}
 $\tilde{J} \leftarrow C^{\sf T}{\sf D}_B^{-1}C$\\
 \For{$k\leftarrow 0$ \KwTo $I_{\rm max}$}{
    $P_r[k] \leftarrow $ from equation \eqref{Eq:Prl-MLDC}\\
    ${\delta}_r[k+1] \leftarrow$ \Solve{${\sf L}_B {\delta}_r = P_r[k]$}\\
      	\uIf{$c = 0$}{
        $\psi[k+1] \leftarrow A_r^{\sf T}{\delta}_r[k+1]$\\  		
        }
    	\uElse{
    	$x \leftarrow$ \Solve{$\tilde{J}(x-x[k]) = C^{\sf T}\barcsin(\psi[k])$}\\
    	$\psi[k+1] \leftarrow A_r^{\sf T}{\delta}_r[k+1] + {\sf D}_{B}^{-1}Cx$\\
  		}
		${\theta}_r[k+1] \leftarrow$ \Solve{$A_r^{\sf T}{\theta}_r = \barcsin({\psi[k+1]})$}\\
  }
  	\KwRet{voltage phase angles ${\theta}_r$}\\
 \caption{Lossy Modified DCPF}
 \label{Alg:LDCPF}
\end{algorithm}

Each iteration for $\psi$ requires solving a linear equation, but with \emph{the same constant coefficient matrix} ${\sf L}_B$, which is sparse and positive-definite. If there are no topology changes, one may compute and store the Cholesky factorization of ${\sf L}_{B}$, and each iteration will require only one forward/backward substitution. Sparsity-exploiting techniques such as optimal bus orderings can also be used for increased speed. Similar statements apply to the linear equation for the $x$-variable updates. The vector $P_{r}[k]$ is given by
\begin{equation}\label{Eq:Prl-MLDC}
\begin{aligned}
P_r[k] &= P_r - G_{\rm diag}V_r^2 + |A|_r\mathsf{D}_{G}\sqrt{\vones[m]-(\psi[k])^2}\,,
\end{aligned}
\end{equation}
and may be thought of as an updated vector of power injections which contains an estimate of the resistive losses, allocated automatically to buses. 
We conclude that for radial systems computing $k$ L-MDCPF iterations has roughly equal complexity to computing $k$ classical DC Power Flows, while for meshed systems an additional $c$-dimensional linear system must be solved to update the $x$ variables.

\smallskip

\subsection{Simulations}
\label{Sec:Simulations}

We now present two simulation studies which illustrate our results. As our primary objective is to compare the Lossy DCPF to the DCPF and to other methods, it is not particularly important whether the models are hot-started with correct voltage magnitudes, or cold-started with $V_i = 1$ for all buses; we have chosen the former.

For the first test, we compare the iterates of the Lossy Modified DCPF \eqref{Eq:MLDCK-PF-General} to the iterates of four other methods: the DCPF \eqref{Eq:DCPF}, the Lossy DCPF \eqref{Eq:LDCK-PF-General}, the Newton-Raphson method (NR), and a chord Newton-Raphson method (CNR); the latter is similar to the fast-decoupled load flow method.\footnote{For a system of nonlinear equations $F(z) = \vzeros[]$, the chord Newton method is defined by the iteration $z_{k+1} = z_{k} - F^\prime(z_0)^{-1}F(z_k)$, where $z_0$ is the initial point for the iteration.} We initialize each method with a flat start $\theta_r = \vzeros[n]$, and measure the error from the exact solution at each iteration using the maximum angle error incurred at any bus $\|\theta_r-\theta_r[k]\|_{\infty}$. 
To tie our simulation results to the theory developed in Section \ref{Sec:LDF-PF3} we consider two networks: the radial IEEE 37 bus system, and the heavily meshed IEEE 118 bus system. Since convergence rates are more easily distinguished in stressed systems, we load both systems along the base case 90\% of the way to power flow insolvability, as determined by continuation power flow.

\begin{figure}[ht!]
\begin{center}
\includegraphics[trim={0.4cm 0.1cm 1.5cm 0.6cm},clip,width = 1\columnwidth]{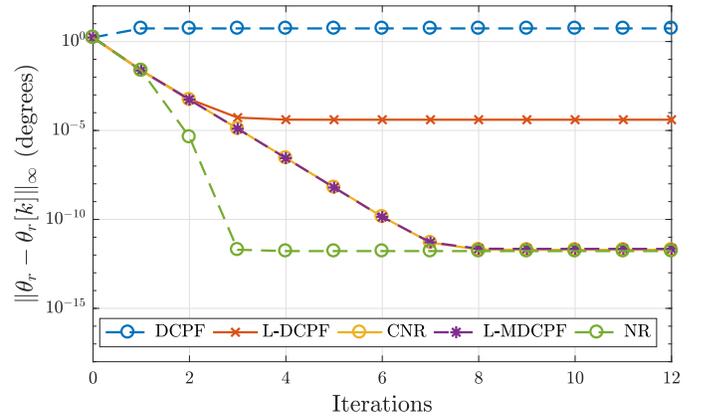}
\caption{Comparison of method iterates for IEEE 37 bus system.}
\label{Fig:Compare37}
\end{center}
\end{figure}

The results for the radial IEEE 37 system are plotted in Figure \ref{Fig:Compare37}. Since the DCPF is non-iterative, it maintains a constant (and rather poor) accuracy over all iterations. The Newton-Raphson iteration shows its characteristic quadratic convergence, taking only a few iterations to reach an accuracy of $10^{-12}$ degrees. The chord Newton-Raphson and the Lossy Modified DCPF plots essentially overlap, taking roughly twice as many iterations as NR. The most interesting plot is that from the Lossy DCPF; the iterates converge to an accuracy of $10^{-4}$ degrees. This is more than accurate enough for practical purposes, but is not the exact solution. For radial networks, the only difference between the L-MDCPF \eqref{Eq:MLDCK-PF-General} and the L-DCPF \eqref{Eq:LDCK-PF-General} is the explicit use of $\barcsin(\cdot)$ in the former. We conclude that the $\barcsin(\cdot)$ nonlinearity provides a mild increase in accuracy, allowing for convergence to the exact solution, but that all methods provide reasonable results in radial systems.

\begin{figure}[ht!]
\begin{center}
\includegraphics[trim={0.4cm 0.1cm 1.5cm 0.62cm},clip,width = 1\columnwidth]{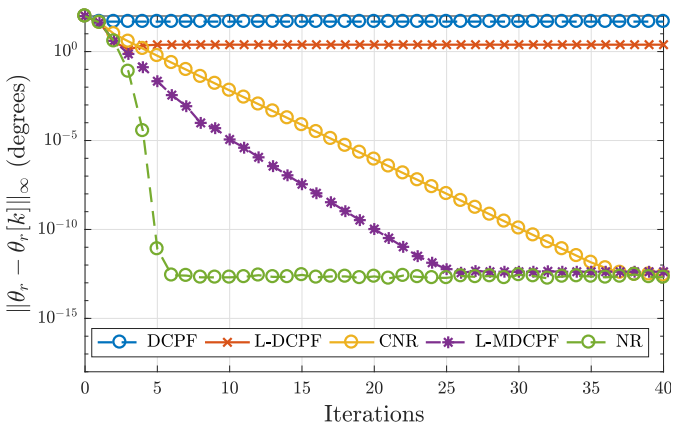}
\caption{Comparison of method iterates for IEEE 118 bus system.}
\label{Fig:Compare118}
\end{center}
\end{figure}

The results for the meshed IEEE 118 system are plotted in Figure \ref{Fig:Compare118}. As before, the DCPF maintains a constant poor accuracy, while the NR iteration converges quadratically. In this example, the Lossy Modified DCPF converges faster than the chord Newton method, taking roughly 10 fewer iterations. Again, the most interesting comparison is between the L-MDCPF and the L-DCPF; the latter converges, but to a point nowhere near the true solution. For meshed networks, the two methods differ significantly due to the extra iteration in the L-MDCPF for the $x$ variables, which enforce Kirchhoff's voltage law. {\tb The trends in Figure \ref{Fig:Compare118} are similar across all meshed test cases, from 39 to 13,659 buses}. We conclude for meshed networks that (i) the L-MDCPF converges linearly to the proper solution, at a rate comparable to fast-decoupled power flow, and that (ii) the L-DCPF iteration does not converge well for meshed networks.

The additional iteration for the $x$-variable in the L-MDCPF \eqref{Eq:MLDCK-PF-General} is required for convergence to the exact solution, but is also undesirable, as it imposes additional computational burden. For our second set of tests, we examine what happens if we simply ignore this part of the iteration, and set $x = \vzeros[c]$ in the L-MDCPF \eqref{Eq:MLDCK-PF-General}. To test over a wider variety of systems, we consider ten standard test cases at base case loading \cite{RDZ-CEM-DG:11,CZ-SF-JM-PP:16}. For each of these networks, Table \ref{Tab:Iterations} displays the maximum phase angle error $\|\theta_r-\theta_r[k]\|_{\infty}$ (in degrees) between the exact solution and the output of the L-MDCPF \eqref{Eq:MLDCK-PF-General} after $k$ iterations. The first column corresponding to $k=1$ is \textemdash{} as we have previously remarked below \eqref{Eq:Psi111} \textemdash{} essentially a version of the MDCPF approximation. Note that some of the MDCPF predictions are extremely inaccurate in large networks. One additional $\psi$-iteration ($k=2$) of the L-MDCPF yields an order-of-magnitude improvement in error, with another iteration ($k=3$) yielding an additional order-of-magnitude in most cases. The accuracy tends to remain constant after the third or fourth iteration; for the PEGASE 13,659 case, $k=4$ brings the error down to 0.5$^\circ$.


\begin{table}[h!]
\begin{center}
\caption{Base Case Testing of Lossy Modified DC Power Flow}
{\renewcommand{\arraystretch}{1}
\begin{tabular}{lccc}
\toprule
Test System & Error (deg) & Error (deg) & Error (deg) \\
& $k=1$ &  $k=2$ & $k=3$ \\
\midrule
New England 39 		& 1.33 & 0.02 & 0.00\\
RTS '96 (2 area) 	& 1.88 & 0.03 & 0.00\\
57 bus system 		& 0.55 & 0.01 & 0.00\\
RTS '96 (3 area)	& 4.16 & 0.06 & 0.01\\
118 bus system		& 3.49 & 0.05 & 0.01\\
300 bus system 		& 19.3 & 0.22 & 0.07\\
Polish 2383wp  		& 5.32 & 0.31 & 0.02\\
PEGASE 2869			& 21.44 & 0.61 & 0.05\\
PEGASE 9241			& 74.05 & 6.02 & 0.37\\
PEGASE 13,659 		& 242.7 & 111.7 & 5.85\\
\bottomrule
\end{tabular}
}
\label{Tab:Iterations}
\end{center}
\end{table}

For heavily loaded networks, we expect that losses will more substantial, and that the gap between the DCPF predictions and the L-DCPF predictions will grow. The previous experiments were repeated, and the results are shown in Table \ref{Tab:Iterations2}.\footnote{PEGASE 13,659 is omitted, as the base case is already more than 90\% of the way to insolvability.} Once again, the first column ($k=1$) shows the wildly inaccurate predictions of the Modified DCPF. One additional iteration of the Lossy DCPF ($k=2$) reduces the error by an order of magnitude, and in some cases by nearly two. A third iteration further improves the results, with the error remaining roughly constant thereafter. We conclude then that even when ignoring the $x$-iteration, the Lossy DCPF is able to correct for resistive losses, and that a result accurate enough for most purposes can be achieved in only one or two iterations beyond the DCPF approximation.
\begin{table}[h!]
\begin{center}
\caption{High Load Testing of Lossy Modified DC Power Flow}
{\renewcommand{\arraystretch}{1}
\begin{tabular}{lccc}
\toprule
Test System & Error (deg) & Error (deg) & Error (deg) \\
& $k=1$ &  $k=2$ & $k=3$ \\
\midrule
New England 39 		& 7.10 & 0.30 & 0.02\\
RTS '96 (2 area) 	& 9.34 & 0.21 & 0.17\\
57 bus system 		& 1.74 & 0.07 & 0.02\\
RTS '96 (3 area)	& 20.5 & 0.42 & 0.2\\
118 bus system		& 42.09 & 3.65 & 1.28\\
300 bus system 		& 35.96 & 0.59 & 0.23\\
Polish 2383wp  		& 19.58 & 2.27 & 0.3\\
PEGASE 2869			& 63.52 & 4.39 & 0.53\\
PEGASE 9241			& 93.37 & 8.07 & 0.51\\
\bottomrule
\end{tabular}
}
\label{Tab:Iterations2}
\end{center}
\end{table}

{\tb Corollary \ref{Cor:RadialConvergence} implies that the L-MDCPF convergence should be insensitive to initialization. As a final test then, we consider the IEEE 118 bus test system with high loading (90\% of the way to insolvability). We generate 1000 random phase angle initial conditions, each with components pulled uniformly from the interval $[-\varphi,\varphi]$ for various values of $\varphi$. For each initial condition, we run the NR, CNR, and L-MDCPF algorithms. If the iterates converge to the known small-angle solution, we mark the test successful; otherwise, we say the solver has failed. Table \ref{Tab:Initial} shows the fraction of successful tests for various values of $\varphi$. For small values of $\varphi$ (i.e., initializations close to $\theta_r = \vzeros[]$) all solvers behave similarly. As $\varphi$ increases, both the NR and CNR solvers increasingly struggle to converge; they either diverge, or converge to an undesirable solution.\footnote{The results are even more extreme in larger test systems, where NR and CNR solvers diverge due to even mild perturbations in initialization.} In contrast, the L-MDCPF iteration recovers the desired solution from every constructed initial condition, even for very highly randomized initial conditions.

\begin{table}[h!]
\begin{center}
\caption{Solver success rates for random initializations (118 bus).}
{\renewcommand{\arraystretch}{1}
\begin{tabular}{|l|c|c|c|}
\toprule
$\varphi$ (deg.) & NR & CNR & L-MDCPF\\
\midrule
15 & 1.00 & 1.00 & 1.00\\
20 & 0.73 & 0.91 & 1.00\\
25 & 0.06 & 0.36 & 1.00\\
30 & 0.00 & 0.02 & 1.00\\
40 & 0.00 & 0.00 & 1.00\\
70 & 0.00 & 0.00 & 1.00\\
80 & 0.00 & 0.00 & 1.00\\
\bottomrule
\end{tabular}
}
\label{Tab:Initial}
\end{center}
\end{table}

}

\section{Conclusions}
\label{Sec:Conclusions}

In this paper we have introduced the Lossy DCPF and the Lossy Modified DCPF, which are generalizations of the classic DC Power Flow and the Modified DC Power Flow to large networks with line resistances. These generalizations yield a convergent sequence of approximations in terms of a sequence of DC Power Flow problems, where the vector of power injections is updated at every step. In standard test cases, the convergence of this sequence is rapid, and taking even one step yields a order-of-magnitude improvement in accuracy over the DC Power Flow. The method is also theoretically justified, with convergence conditions given for radial networks.

There are many avenues for future work. Theoretically, it would be desirable (i) to obtain less conservative conditions in Theorem \ref{Thm:Conditions}, (ii) to obtain convergence conditions for \eqref{Eq:MLDCK-PF-General} in meshed networks, and (iii) to further extend the algorithm and theoretical results from active power flow to full, coupled AC power flow. Further analysis and understanding of iterative solution methods for meshed networks is needed. On the applications side, the Lossy DCPF algorithm seems well-suited for market and contingency analysis, and future work will examine whether the L-DCPF can enable rapid yet accurate calculation of locational marginal prices \cite{FL-RB:07} and line outage distribution factors, especially in large, heavily loaded networks.

\appendices

\section{Technical Proofs}
\label{App:1}

\begin{pfof}{Proposition \ref{Thm:GeneralConvergence}}
If $\psi[k] \rightarrow \psi^*$ and $x[k] \rightarrow x^*$, then by construction $(\psi^*,x^*)$ together satisfy \eqref{Eq:PsiFP-General}. Substitution using \eqref{Eq:PsiFP-General} shows immediately that  $\psi^*$ solves \eqref{Eq:ActivePsia}. Note that, by necessity, $\psi^*$ must satisfy $\|\psi^*\|_{\infty} < 1$, since otherwise the right-hand side of \eqref{Eq:PsiFP-General} would be complex while the left-hand side of \eqref{Eq:PsiFP-General} would be real. By construction, $\psi^*$ also satisfies \eqref{Eq:PsiFP-General2}, or equivalently $C^{\sf T}\barcsin(\psi^*) = 2\pi \xi$ for some integer vector $\xi = (\xi_1,\ldots,\xi_c)^{\sf T}$. Since $C$ has full rank and is totally unimodular \cite[Theorem 3.4]{TK-CL-KM-DM-RR-TU-KAZ:09}, we can always find another integer vector $\xi^\prime$ such that $C^{\sf T}\xi^\prime = \xi$ \cite[Theorem 5.20]{BK-JV:08}.
The general solution to $C^{\sf T}\barcsin(\psi^*) = 2\pi \xi$ may therefore be written as $\barcsin(\psi^*) = 2\pi \xi^\prime + A_r^{\sf T}\theta_r$ for some $\theta_r \in \real^{n}$, with $2\pi \xi^\prime$ the particular solution and $A_r^{\sf T}\theta_r$ parameterizing the homogeneous solution. Taking the $\bsin(\cdot)$ of both sides shows that $\psi^*$ satisfies \eqref{Eq:ActivePsib}, and the conclusion follows.\end{pfof}

\begin{pfof}{Theorem \ref{Thm:Conditions}}
Beginning with the auxiliary equation \eqref{Eq:ActivePsia}, we add and subtract $|A|_r{\sf D}_{G}\vones[m]$ to obtain
\begin{equation}\label{Eq:ActivePsi-App1}
\begin{aligned}
P_r &= A_r\mathsf{D}_{B}\psi + |A|_r\mathsf{D}_{G}\left(\vones[m]-\sqrt{\vones[m]-\psi^2}\right)\\ & \qquad+ G_{\rm diag}V_r^2 - |A|_r\mathsf{D}_{G}\vones[m]
\,.
\end{aligned}
\end{equation}
In components, the final two terms read as
\begin{align*}
G_{ii}V_i^2 - \sum_{j=1}^{n+1}\nolimits V_iV_jG_{ij} = -\sum_{j\neq i}^n\nolimits V_iG_{ij}(V_i-V_j)\,,
\end{align*}
where we inserted $G_{ii} = -\sum_{j=1}^{n+1}G_{ij}$ since there are no off-nominal tap ratios. Since by assumption $V_i = V \in \real$ for all $i \in \{1,\ldots,n+1\}$, this expression is identically zero, and these terms vanish from \eqref{Eq:ActivePsi-App1}. 
Since the network is radial, $\mathrm{ker}(A_r) = \emptyset$ \cite[Prop. 4.3]{NB:94}. It follows that $A_r$ is invertible, and we may rearrange \eqref{Eq:ActivePsi-App1} to obtain
\begin{equation}\label{Eq:FPEqa}
\psi = {\sf D}_{B}^{-1}A_r^{-1}\left(P_r + |A|_r\mathsf{D}_{G}\left(\vones[m]-\sqrt{\vones[m]-\psi^2}\right)\right)\,.
\end{equation}
We now claim that ${\sf D}_B^{-1}A_r^{-1} = A_r^{\sf T}{\sf L}_{B}^{-1}$. To see this, let $X = A_r^{\sf T}{\sf L}_{B}^{-1}$. Then $A_r{\sf D}_{B}X = (A_r{\sf D}_BA_r^{\sf T}){\sf L}_B^{-1} = {\sf L}_{B}{\sf L}_B^{-1} = I_{m}$, and the result follows. Therefore, \eqref{Eq:FPEqa} may be written as
\begin{align}
\psi &= \underbrace{A_r^{\sf T}{\sf L}_B^{-1}P_r}_{\triangleq \psi_{\rm MDC}} + \underbrace{{\sf D}_{B}^{-1}A_r^{-1}|A|_r\mathsf{D}_{G}}_{\triangleq H}\left(\vones[m]-\sqrt{\vones[m]-[\psi]\psi}\right)\nonumber\\
&= \psi_{\rm MDC} + H\left(\vones[m]-\sqrt{\vones[m]-\psi^2}\right) \triangleq f(\psi)\,.
\label{Eq:FPEquation}
\end{align}
In other words, $\psi$ solves the auxiliary equation \eqref{Eq:ActivePsia} if and only if $\psi$ solves \eqref{Eq:FPEquation}, which is a fixed point equation $\psi = f(\psi)$ for the variable $\psi$. We will first seek to find a closed invariant set for $f$, i.e., a closed set $\mathcal{I}$ such that $\map{f}{\mathcal{I}}{\mathcal{I}}$. For $\beta \in [0,1)$, let $\mathcal{I}(\beta) = \setdef{\psi \in \real^m}{\|\psi\|_{\infty} \leq \beta}$, and suppose that $\psi \in \mathcal{I}(\beta)$. Then
\begin{align}
\|f(\psi)\|_{\infty} &\leq \|\psi_{\rm MDC}\|_{\infty} + \|H\|_{\infty}(1-\sqrt{1-\beta^2})\nonumber\\
\label{Eq:fbound}
&\leq \Gamma + \rho(1-\sqrt{1-\beta^2})\,,
\end{align}
where we have inserted the definition of $\rho$ from \eqref{Eq:rho} and $\Gamma$ from \eqref{Eq:Gamma}. We now require that \eqref{Eq:fbound} is further upper bounded by $\beta$, which will imply that $f(\psi) \in \mathcal{I}(\beta)$. This leads to the inequality
$$
\Gamma + \rho(1-\sqrt{1-\beta^2}) \leq \beta\,,
$$
which can be easily manipulated into
\begin{equation}\label{Eq:ScalarInequality}
(1+\rho^2)\beta^2 - 2(\Gamma+\rho)\beta + (\Gamma+\rho)^2-\rho^2 \leq 0\,.
\end{equation}
While we omit the details, one can solve \eqref{Eq:ScalarInequality} with \emph{equality} sign for $\beta$, obtaining the two solutions
\begin{equation}\label{Eq:Betapm}
\beta_{\pm} = \frac{\Gamma + \rho}{1+\rho^2} \pm \frac{\rho}{1+\rho^2}\sqrt{1-(\Gamma^2+2\Gamma\rho)}\,.
\end{equation}
The solutions are well-defined and distinct if and only if the stated condition \eqref{Eq:Condition} holds, in which case they satisfy $0 \leq \beta_- < \beta_+ \leq 1$. 
\begin{figure}[t!]
\begin{center}
\includegraphics[width = 0.5\columnwidth]{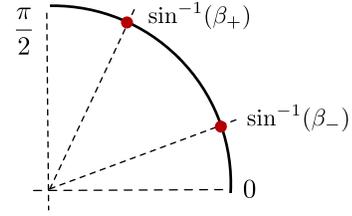}
\caption{Two points on the circle defined by $\beta_-$ and $\beta_+$.}
\label{Fig:beta}
\end{center}
\end{figure}
Inspection then reveals that \eqref{Eq:ScalarInequality} is satisfied for all $\beta \in [\beta_-, \beta_+]$, and satisfied with strict inequality sign for $\beta \in (\beta_-,\beta_+)$. In the latter case, pick any $\psi \in \mathcal{I}(\beta)$ such that $\|\psi\|_{\infty} = \beta$. Then $\|f(\psi)\|_{\infty} < \beta = \|\psi\|_{\infty}$, and therefore $f(\psi) \neq \psi$. It follows that there can exist no fixed points $\psi$ with components $\psi_i \in (\beta_-,\beta_+)$, and therefore there exist no power flow solutions in the set \eqref{Eq:NoSolutions}. In particular, we have shown that if $\psi \in \mathcal{I}(\beta_-)$, then $f(\psi) \in \mathcal{I}(\beta_-)$, so $\mathcal{I}(\beta_-)$ is a compact invariant set for $f$. We now show that $f$ is a contraction on $\mathcal{I}(\beta_-)$. From \eqref{Eq:FPEquation}, we have
\begin{equation}
\frac{\partial f}{\partial \psi}(\psi) = H(I_m-[\psi]^2)^{-\frac{1}{2}}[\psi]\,.
\end{equation}
If $\psi \in \mathcal{I}(\beta_-)$, then it follows that
\begin{equation}\label{Eq:AbsJ}
\Big\| \frac{\partial f}{\partial \psi}(\psi)\Big\|_{\infty} \leq \|H\|_{\infty}\frac{\beta_-}{\left(1-\beta_-^2\right)^{1/2}} = \underbrace{\rho \frac{\beta_-}{\left(1-\beta_-^2\right)^{1/2}}}_{\triangleq c(\rho,\Gamma)}\,.
\end{equation}
We now seek to show that $c(\rho,\Gamma) < 1$ for all $(\Gamma,\rho)$ satisfying the condition \eqref{Eq:Condition}. First note that $c(0,0) = 0$. Next note that since $\Gamma \mapsto \beta_-$ is strictly increasing, so is $\Gamma \mapsto c(\rho,\Gamma)$. While less obvious, it is also true that $\rho \mapsto \beta_-$ is strictly increasing, and therefore so is $\rho \mapsto c$. It therefore suffices to check that $c(\rho,\Gamma) \leq 1$ when the condition \eqref{Eq:Condition} is satisfied with equality sign. Solving $\Gamma^2 + 2\Gamma\rho - 1 = 0$ for $\Gamma > 0$ yields $\Gamma_{\rm crit} = -\rho + \sqrt{1+\rho^2}$. Some straightforward algebra then shows that $c(\rho,\Gamma_{\rm crit}) = 1$, and therefore $\|\frac{\partial f}{\partial \psi}(\psi)\|_{\infty} < 1$ for all $(\Gamma,\rho)$ satisfying the condition \eqref{Eq:Condition}. 
Since $\mathcal{I}(\beta_-)$ is convex, it follows by standard results that $f$ is a contraction mapping on $\mathcal{I}(\beta_-)$. Therefore, the Banach Fixed-Point Theorem \cite[Theorem 9.32]{WR:76} guarantees that $f$ possess a unique fixed point $\psi^*$ satisfying $\|\psi^*\|_{\infty} \leq \beta_- < 1$. Therefore, $\psi^*$ is the unique solution contained in the set $\mathcal{I}(\beta_-)$ to the auxiliary equation \eqref{Eq:ActivePsia}. Since $\|\psi^*\|_{\infty} < 1$, we may apply $\arcsin(\cdot)$ component-wise to \eqref{Eq:ActivePsib} to obtain
\begin{equation}\label{Eq:InvertedAngle}
A_r^{\sf T}\theta_r^* = \barcsin(\psi^*).
\end{equation}
Since the network is radial, $\mathrm{ker}(A_r) = \emptyset$, and therefore $\mathrm{im}(A_r^{\sf T}) = \real^{n}$ \cite[Prop. 4.3]{NB:94}. The right-hand side of \eqref{Eq:InvertedAngle} is therefore always in the image of $A_r^{\sf T}$, and we may uniquely solve \eqref{Eq:InvertedAngle} for $\theta_r^* = (A_rA_r^{\sf T})^{-1}A_r\barcsin(\psi^*)$.
\end{pfof}

\IEEEpeerreviewmaketitle


\ifCLASSOPTIONcaptionsoff
  \newpage
\fi


\bibliographystyle{IEEEtran}
\bibliography{alias,JWSP,New,Main}

\begin{IEEEbiography}[{\includegraphics[width=1in,height=1.25in,clip,keepaspectratio]{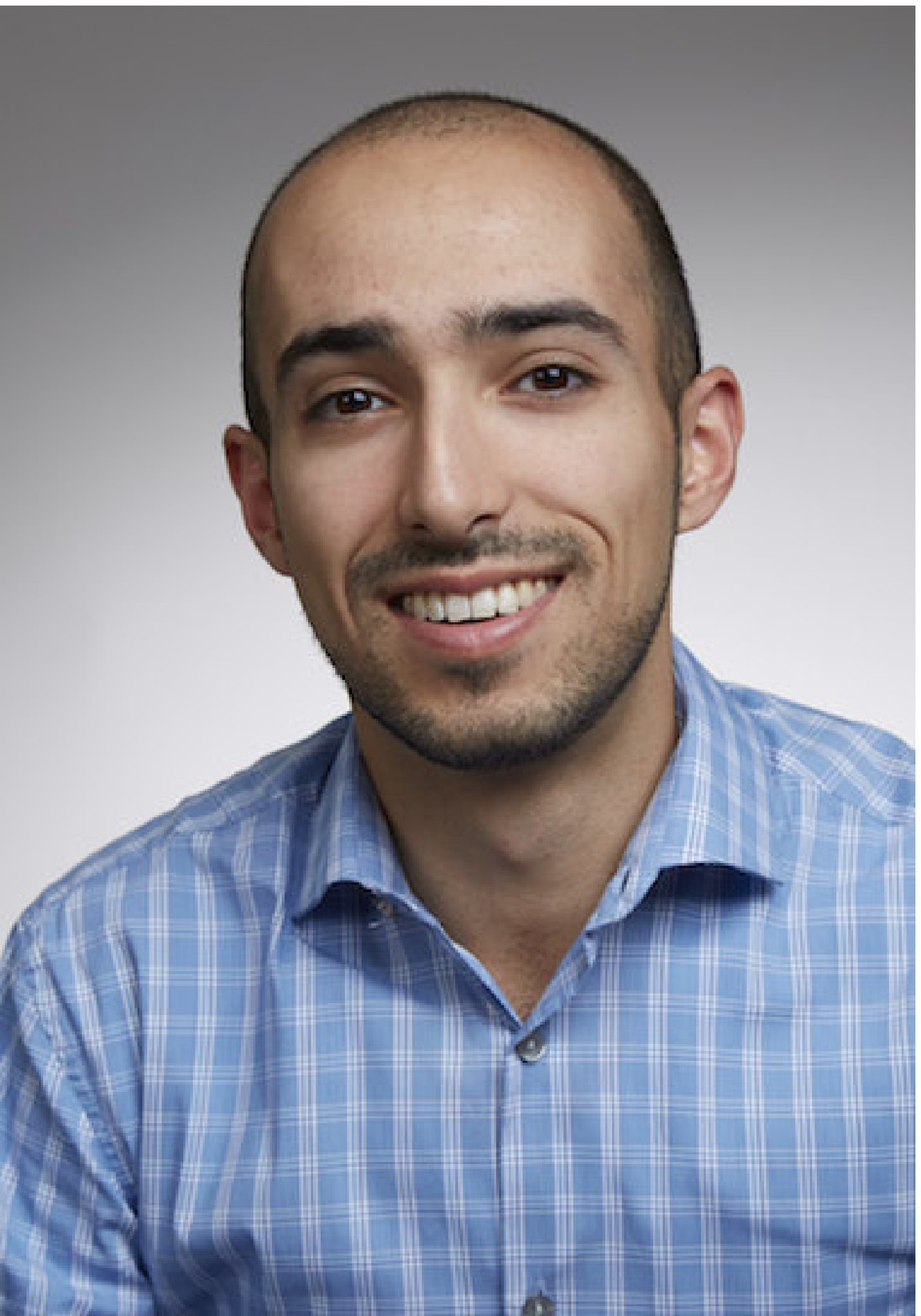}}]{John W. Simpson-Porco} (S'11--M'16) received the B.Sc. degree in engineering physics from Queen's University, Kingston, ON, Canada in 2010, and the Ph.D. degree in mechanical engineering from the University of California at Santa Barbara, Santa Barbara, CA, USA in 2015.

He is currently an Assistant Professor of Electrical and Computer Engineering at the University of Waterloo, Waterloo, ON, Canada. He was previously a visiting scientist with the Automatic Control Laboratory at ETH Z\"{u}rich, Z\"{u}rich, Switzerland. His research focuses on the control and optimization of multi-agent systems and networks, with applications in modernized power grids.

Prof. Simpson-Porco is a recipient of the 2012--2014 IFAC Automatica Prize and the Center for Control, Dynamical Systems and Computation Best Thesis Award and Outstanding Scholar Fellowship.
\end{IEEEbiography}

\end{document}